\documentclass{article}
\newtheorem{theorem}{Theorem}

\begin{document}

\title{Partitions with parts in a finite set
\thanks{
Supported in part by grants from the PSC--CUNY Research Award Program and the
NSA Mathematical Sciences Program.}}
\author{Melvyn B. Nathanson\\
Department of Mathematics\\
Lehman College (CUNY)\\
Bronx, New York 10468\\
e-mail: nathansn@alpha.lehman.cuny.edu}
\date{}

\maketitle
\begin{abstract}
Let $A$ be a nonempty finite set of relatively prime positive integers,
and let $p_A(n)$ denote the number of partitions of $n$ with parts in $A$.
An elementary arithmetic argument is used to prove
the asymptotic formula
\[
p_A(n) = \left(\frac{1}{\prod_{a\in A}a}\right) \frac{n^{k-1}}{(k-1)!}
+ O\left( n^{k-2}\right).                                  
\]
\end{abstract}

Let $A$ be a nonempty set of positive integers.
A {\em partition} of a positive integer $n$
with parts in $A$ is a representation of $n$
as a sum of not necessarily distinct elements of $A$.
Two partitions are considered the same if they differ
only in the order of their summands.
The {\em partition function} of the set $A$,
denoted $p_A(n)$, counts the number of partitions of $n$ with parts in $A$.

If $A$ is a finite set of positive integers with no common factor
greater than 1, then every sufficiently large integer can be written
as a sum of elements of $A$ (see Nathanson~\cite{nath72f}
and Han, Kirfel, and  Nathanson~\cite{han-kirf-nath98}), and so
$p_A(n) \geq 1$ for all $n \geq n_0.$
In the special case that $A$ is the set of the first $k$ integers,
it is known that
\[
p_A(n) \sim \frac{n^{k-1}}{k!(k-1)!}.
\]
Erd\H os and Lehner\cite{erdo-lehn41} proved that this asymptotic formula
holds uniformly for $k=o(n^{1/3})$.
If $A$ is an arbitrary finite set of relatively prime positive integers,
then
\begin{equation}   \label{kparts:formula}
p_A(n) \sim \left(\frac{1}{\prod_{a\in A}a}\right) \frac{n^{k-1}}{(k-1)!}.
\end{equation}
The usual proof of this result
(Netto~\cite{nett27}, P\' olya--Szeg\" o~\cite[Problem 27]{poly-szeg25})
is based on the partial fraction
decomposition of the generating function for $p_A(n)$.
The purpose of this note is to give a simple, purely arithmetic proof
of~(\ref{kparts:formula}).

We define $p_A(0) = 1.$

\begin{theorem}                    \label{para:Aset}
Let $A = \{a_1,\ldots, a_k\}$ be a set of $k$ relatively prime
positive integers, that is,
\[
\gcd(A) = (a_1,\ldots,a_k) = 1.
\]
Let $p_A(n)$ denote the number of partitions of $n$ into parts
belonging to $A$.
Then
\[
p_A(n) = \left(\frac{1}{\prod_{a\in A}a}\right) \frac{n^{k-1}}{(k-1)!}
+ O\left( n^{k-2}\right).
\]
\end{theorem}

{\bf Proof.}
Let $k = |A|.$
The proof is by induction on $k$.
If $k = 1$, then $A = \{1\}$
and
\[
p_A(n) = 1,
\]
since every positive integer has a unique partition into a sum of 1's.

Let $k \geq 2,$ and assume that the Theorem holds for $k-1$.
Let
\[
d = (a_1,\ldots, a_{k-1}).
\]
Then
\[
(d,a_k) = 1.
\]
For $i = 1,\ldots, k-1$, we set
\[
a_i' = \frac{a_i}{d}.
\]
Then
\[
A' = \{a'_1,\ldots, a'_{k-1}\}
\]
is a set of $k-1$ relatively prime positive integers,
that is,
\[
\gcd(A') = 1.
\]

Since the induction assumption holds for $A'$, we have
\[
p_{A'}(n) = \left( \frac{1}{\prod_{i=1}^{k-1}a'_i}\right) \frac{n^{k-2}}{(k-2)!}
+ O\left( n^{k-3} \right)
\]
for all nonnegative integers $n.$

Let $n \geq (d-1)a_k$.  Since $(d,a_k) = 1$, there exists a unique integer
$u$ such that $0 \leq u \leq d-1$ and
\[
n \equiv ua_k\pmod{d}.
\]
Then
\[
m = \frac{n-ua_k}{d}
\]
is a nonnegative integer, and
\[
m = O(n).
\]
If $v$ is any nonnegative integer such that
\[
n \equiv va_k\pmod{d},
\]
then $va_k \equiv ua_k\pmod{d}$ and so $v \equiv u \pmod{d}$,
that is, $v = u+\ell d$ for some nonnegative integer $\ell.$
If
\[
n-va_k = n - (u+\ell d)a_k \geq 0,
\]
then
\[
0 \leq \ell \leq \left[ \frac{n}{da_k} - \frac{u}{d} \right]
= \left[\frac{m}{a_k}\right] = r.
\]
We note that
\[
r = O(n).
\]
Let $\pi$ be a partition of $n$ into parts belonging to $A.$
If $\pi$ contains exactly $v$ parts equal to $a_k$,
then $n-va_k \geq 0$ and $n - va_k\equiv 0\pmod{d}$, since $n-va_k$
is a sum of elements in $\{a_1, \ldots, a_{k-1}\}$,
and each of the elements in this set is divisible by $d$.
Therefore, $v = u+\ell d$, where $0 \leq \ell \leq r.$
Consequently, we can divide the partitions of $n$ with parts
in $A$ into $r+1$ classes, where, for each $\ell = 0,1,\ldots, r,$
a partition belongs to class $\ell$ if it contain exactly
$u + \ell d$ parts equal to $a_k$.
The number of partitions of $n$ with exactly $u+\ell d$ parts equal to $a_k$
is exactly the number of partitions of $n - (u+\ell d)a_k$ into parts
belonging to the set $\{a_1,\ldots, a_{k-1}\}$, or, equivalently,
the number of partitions of
\[
\frac{n - (u+\ell d)a_k}{d}
\]
into parts belonging to $A'$, which is exactly
\[
p_{A'}\left( \frac{n - (u+\ell d)a_k}{d} \right)
= p_{A'}\left( m - \ell a_k \right).
\]
Therefore,
\begin{eqnarray*}
p_A(n)
& = & \sum_{\ell = 0}^r p_{A'}(m-\ell a_k)  \\
& = & \left(\frac{1}{\prod_{i=1}^{k-1}a'_i}\right) \sum_{\ell = 0}^{r}
\left(\frac{ (m-\ell a_k)^{k-2}}{(k-2)!} + O(m^{k-3}) \right)  \\
& = & \left(\frac{d^{k-1}}{\prod_{i=1}^{k-1}a_i}\right)
\sum_{\ell = 0}^{r} \frac{ (m-\ell a_k)^{k-2}}{(k-2)!} + O(n^{k-2}).
\end{eqnarray*}
To evaluate the inner sum, we note that
\[
\sum_{\ell = 0}^{r} \ell^j
= \frac{r^{j+1}}{(j+1)} + O(r^j)
\]
and
\[
\sum_{j=0}^{k-2} (-1)^j {k-1\choose j+1}
= - \sum_{j=1}^{k-1} (-1)^j {k-1\choose j} = 1.
\]
Then
\begin{eqnarray*}
\sum_{\ell = 0}^{r} \frac{ (m-\ell a_k)^{k-2}}{(k-2)!}
& = & \frac{1}{(k-2)!} \sum_{\ell = 0}^{r}
\sum_{j=0}^{k-2} {k-2\choose j} m^{k-2-j}(-\ell a_k)^j  \\
& = & \frac{1}{(k-2)!} \sum_{j=0}^{k-2}{k-2\choose j} m^{k-2-j}(-a_k)^j
\sum_{\ell = 0}^{r} \ell^j \\
& = & \frac{1}{(k-2)!} \sum_{j=0}^{k-2}{k-2\choose j} m^{k-2-j}(-a_k)^j
\left( \frac{r^{j+1}}{(j+1)} + O(r^j)\right)  \\
& = & \frac{1}{(k-2)!} \sum_{j=0}^{k-2}{k-2\choose j} m^{k-2-j}(-a_k)^j
\left( \frac{m^{j+1}}{a_k^{j+1}(j+1)} + O(m^j)\right)  \\
& = & \frac{m^{k-1}}{a_k}
\sum_{j=0}^{k-2}{k-2\choose j}\frac{(-1)^j}{(k-2)!(j+1)} + O(m^{k-2}) \\
& = & \frac{m^{k-1}}{a_k}
\sum_{j=0}^{k-2}\frac{(-1)^j}{(k-2-j)!j!(j+1)} + O(m^{k-2}) \\
& = & \frac{m^{k-1}}{a_k}
\sum_{j=0}^{k-2}\frac{(-1)^j}{(k-1-(j+1))!(j+1)!} + O(m^{k-2}) \\
& = & \frac{m^{k-1}}{a_k(k-1)!} \sum_{j=0}^{k-2} (-1)^j {k-1\choose j+1}
+ O(m^{k-2}) \\
& = & \frac{m^{k-1}}{a_k(k-1)!} + O(m^{k-2}).
\end{eqnarray*}
Therefore,
\begin{eqnarray*}
p_A(n)
& = & \left(\frac{d^{k-1}}{\prod_{i=1}^{k-1}a_i}\right)
\sum_{\ell = 0}^{r} \frac{ (m-\ell a_k)^{k-2}}{(k-2)!} + O(n^{k-2})  \\
& = & \left(\frac{d^{k-1}}{\prod_{i=1}^{k-1}a_i}\right)
\left( \frac{m^{k-1}}{a_k(k-1)!} + O(n^{k-2}) \right) + O(n^{k-2}) \\
& = & \left(\frac{d^{k-1}}{\prod_{i=1}^{k-1}a_i}\right)
\left( \frac{1}{a_k(k-1)!}\right)
\left( \frac{n}{d}-\frac{ua_k}{d} \right)^{k-1} + O(n^{k-2}) \\
& = & \left(\frac{d^{k-1}}{\prod_{i=1}^{k-1}a_i}\right)
\left( \frac{1}{a_k(k-1)!}\right)\left( \frac{n}{d}\right)^{k-1} + O(n^{k-2}) \\
& = & \left(\frac{1}{\prod_{i=1}^{k}a_i}\right)
\frac{n^{k-1}}{(k-1)!} + O(n^{k-2}).
\end{eqnarray*}
This completes the proof.

\end{document}